\documentclass[a4paper,11pt,french,epic,eepic]{article}
\usepackage[latin1]{inputenc}
\usepackage{amsmath,amssymb,amsthm}
\usepackage[dvips]{graphicx,epsfig}
\newcommand{\GG}{\mathbb{G}}
\newcommand{\NN}[0]{\mathbb{N}}
\newcommand{\PP}{I\!\!P}

\newcommand{\ZZ}{\mbox{\rm \lower0.3pt\hbox{$\angle\!\!\!$}Z}}

\newtheorem{nt}{Notation}
\newtheorem{conjecture}[nt]{Conjecture}
\newtheorem{coro}[nt]{Corollaire}
\newtheorem{defi}[nt]{D\'efinition}

\newtheorem{ex}[nt]{Exemple}
\newtheorem{lm}[nt]{Lemme}
\newtheorem{prop}[nt]{Proposition}
\newtheorem{rappel}[nt]{Rappel}
\newtheorem{rem}[nt]{Remarque}
\newtheorem{thm}[nt]{Th\'eor\`eme}

\newcommand{\f}{\frac}

\newcommand{\fd}{\ensuremath{\rightarrow}}

\newcommand{\findem}{\begin{flushright}\rule{2mm}{2mm}
\end{flushright}}

\newcommand{\fxc}{\ensuremath{\mathcal C}}

\renewcommand{\phi}{\ensuremath{\varphi}}

\newcommand{\inc}{\ensuremath{\subset}}
\newcommand{\nl}{\ \\[2mm]}

\newcommand{\ox}{\otimes }

\newcommand{\plp}{\PP^2}
\renewcommand{\phi}{\ensuremath{\varphi}}

\newcommand{\s}{Spec\;}

\newcommand{\x}{\ensuremath{\times}}

\newcommand{\barr}{\begin{array}}
\newcommand{\earr}{\end{array}}
\newcommand{\bit}{\begin{itemize}}
\newcommand{\eit}{\end{itemize}}
\newcommand{\beq}{\begin{eqnarray*}}
\newcommand{\eeq}{\end{eqnarray*}}
\newcommand{\beqn}{\begin{eqnarray}}
\newcommand{\eeqn}{\end{eqnarray}}
\newcommand{\bconj}{\begin{conjecture}}
\newcommand{\econj}{\end{conjecture}}
\newcommand{\bcor}{\begin{coro}}
\newcommand{\ecor}{\end{coro} \noindent}
\newcommand{\ben}{\begin{enumerate}}
\newcommand{\een}{\end{enumerate} \noindent}
\newcommand{\bnot}{\begin{nt} }
\newcommand{\enot}{\end{nt} \noindent}
\newcommand{\bdefi}{\begin{defi}}
\newcommand{\edefi}{\end{defi} \noindent}
\newcommand{\bprop}{\begin{prop}}
\newcommand{\eprop}{\end{prop} \noindent \textit{D\'emonstration: }}
\newcommand{\brap}{\begin{rappel}}
\newcommand{\erap}{\end{rappel} \noindent }
\newcommand{\brq}{\begin{rem}}
\newcommand{\erq}{\end{rem} \noindent }
\newcommand{\bthm}{\begin{thm}}
\newcommand{\ethm}{\end{thm}\noindent \textit{D\'emonstration: }}
\newcommand{\blm}{\begin{lm}}
\newcommand{\elm}{\end{lm} \noindent }
\newcommand{\bex}{\begin{ex}}
\newcommand{\eex}{\end{ex}\noindent }

\font \tengothic=eufm10 scaled\magstep 1
\font\sevengothic=eufm7 scaled\magstep 1
\newfam\gothicfam
      \textfont\gothicfam=\tengothic
      \scriptfont\gothicfam=\sevengothic
\def\goth#1{{\fam\gothicfam #1}}

\newcommand{\fxl}{{\cal L}}

\title{Dimension des syst\`emes lin\'eaires: une approche diff\'erentielle
et combinatoire}
\author{Laurent Evain}
\date{evain@tonton.univ-angers.fr}

\begin{document}
\maketitle \noindent
{\bf R\'esum\'e}: On d\'emontre un th\'eor\`eme de majoration pour 
la dimension de certains syst\`emes lin\'eaires.
Ce th\'eor\`eme am\'eliore la m\'ethode d'Horace diff\'e\-ren\-ti\-elle 
introduite par Alexander-Hirschowitz, et avait \'et\'e  conjectur\'e 
par Simpson.
Les applications envisag\'ees sont  le calcul de dimension
de syst\`emes lin\'eaires d'hypersurfaces de $\PP^n$
\`a singularit\'es g\'en\'eriques impos\'ees et le calcul de collisions 
de gros points dans $\plp$. Ces applications 
seront trait\'ees ind\'ependamment de ce papier, mais un exemple 
simple trait\'e dans l'intro\-duc\-tion laisse deviner comment le 
th\'eor\`eme sera utilis\'e.

\section{Introduction par un exemple}
Consid\'erons le syst\`eme lin\'eaire $\fxl_t$ des courbes projectives
planes de degr\'e $d$ passant par trois points fixes   
$p_1,p_2,p_3$ et par un point $p_4(t)$ avec multiplicit\'es respectives 
$m_1,m_2,m_3$ et $m_4$. Supposons que  
$p_1,p_2,p_3$ soient align\'es sur une droite, et  que 
$p_4(t)$ soit g\'en\'erique dans le plan. 
Le syst\`eme est de dimension projective au moins $\f{d(d+3)}{2} 
-\sum \f{m_i(m_i+1)}{2}$. 
Et la dimension est exactement
$\f{d(d+3)}{2} 
-\sum \f{m_i(m_i+1)}{2}$
si les conditions impos\'ees par les points 
multiples sont ind\'ependantes. 
\nl
Choisissons dans notre exemple les conditions num\'eriques 
$m_1=m_2=m_3=1$, $m_4=3$ et $d=5$. 
On veut montrer que le syst\`eme $\fxl_t$
est de dimension onze, et il suffit de voir qu'il est de dimension 
au plus onze. On sp\'ecialise le point g\'en\'erique $p_4(t)$ en 
un point $p_4(0)$ de la droite  $D$ joignant 
$p_1,p_2$ et $p_3$, ce qui  d\'efinit un syst\`eme lin\'eaire $\fxl_0$.
Par semi-continuit\'e, $dim \fxl_t \leq dim \fxl_0$.  
Les diviseurs de $\fxl_0$ sont des courbes de degr\'e cinq qui coupent
la droite $D$ le long d'un sch\'ema ponctuel 
de degr\'e six, donc ils contiennent 
$D$. En soustrayant $D$ \`a chaque diviseur de $\fxl_0$, 
on voit que la dimension de $\fxl_0$  est la m\^eme que celle du 
syst\`eme lin\'eaire des courbes de degr\'e quatre passant par $p_4(0)$ 
avec multiplicit\'e deux, c'est \`a dire onze. On avait donc bien 
$dim \fxl_t=11$.
\nl
Il existe de nombreuses situations  dans lesquelles
on essaie
d'appliquer la strat\'egie  pr\'ec\'edente: on sp\'ecialise 
des points g\'en\'eriques sur des diviseurs de sorte que le probl\`eme 
se simplifie en position sp\'eciale et on conclut par un 
argument de semi-continuit\'e.
Bien s\^ur les conditions num\'eriques de l'exemple ont \'et\'e choisies 
pour  que la strat\'egie s'applique  sans difficult\'e. En revanche, 
il existe en g\'en\'eral des difficult\'es  num\'eriques, comme 
l'illustre le cas
suivant.
\nl
Choisissons dans  notre  exemple introductif 
$m_1=m_2=m_3=2$, $m_4=3$, et $d=6$.
On veut montrer que le syst\`eme $\fxl_t$ est de dimension douze.
Sp\'ecialisons  le point $p_4(t)$ en un point $p_4(0)$ de 
la droite $D$ joignant $p_1,p_2$ et 
$p_3$. Comme pr\'ec\'edemment, $D$ est une composante du 
syst\`eme lin\'eaire $\fxl_0$ dans cette position sp\'eciale,
donc $\fxl_0$ 
a la m\^eme dimension que le syst\`eme des courbes de degr\'e 
cinq passant par $p_1,p_2,p_3,p_4(0)$ avec multiplicit\'e 
un, un, un et deux, c'est \`a dire au moins quatorze. On ne peut pas conclure. 
En fait, la dimension du syst\`eme lin\'eaire en position 
sp\'eciale a saut\'e car on a mis ``trop de conditions sur la droite'':
il suffit qu'un diviseur $\Delta$  de degr\'e six coupe la droite $D$
le long d'un sch\'ema de degr\'e sept pour que $D$ soit inclus dans $\Delta$, 
or un diviseur de $\fxl_0$ coupe la droite
le long d'un sch\'ema de degr\'e neuf.
\nl
La m\'ethode d'Horace
[AH1,AH2,H] propose un ensemble de techniques 
pour g\'erer les probl\`emes num\'eriques qui apparaissent 
lorsqu'on traite des exemples pr\'ecis. On se propose dans ce papier  
d'enrichir la m\'ethode d'Horace d'un nouveau th\'eor\`eme (Th\'eor\`eme
\ref{theo}).
\nl
Alors que l'\'enonc\'e g\'en\'eral n\'ecessite quelques notations, on 
peut illustrer facilement le th\'eor\`eme sur l'exemple pr\'ec\'edent.
\nl
Quand $p_4(t)$ n'est pas 
sur $D$, un diviseur $\Delta$ du syst\`eme lin\'eaire $\fxl_t$
coupe $D$ deux fois en $p_1$, deux fois en $p_2$ et deux fois 
en $p_3$. 
Il manque encore une condition sur la droite pour que $D$ soit 
composante fixe de $\fxl_t$. Raisonnons malgr\'e tout comme si  
$D$ \'etait  composante fixe. Alors, un diviseur de $\fxl_t(-D)$ 
serait 
une courbe de degr\'e 
cinq qui couperait $D$ en $p_1,p_2,p_3$. Il manquerait cette fois-ci 
trois conditions pour que $D$ soit composante du syst\`eme $\fxl_t(-D)$,
c'est \`a dire pour que $2D$ soit dans le lieu fixe de $\fxl$.
\nl
On  va ``prendre les conditions dont on a besoin 
sur le point $p_4$'', qui est un point multiple d'ordre  trois,
lorsque celui-ci approche de $D$.
On pr\'el\`eve les conditions par l'op\'eration 
combinatoire suivante. Passer par 
un point de multiplicit\'e trois \'equivaut \`a contenir un 
gros point de  taille trois de $\plp$. Un gros point de taille 
trois est un sch\'ema mon\^omial i.e.
d\'efini par des \'equations mon\^omiales dans un bon syst\`eme 
de coordonn\'ees. On associe des objets combinatoires
aux sch\'emas 
mon\^omiaux: des escaliers.  
Dans le cas du gros point de taille trois,  l'escalier associ\'e est 
dessin\'e dans la figure ci-apr\`es.
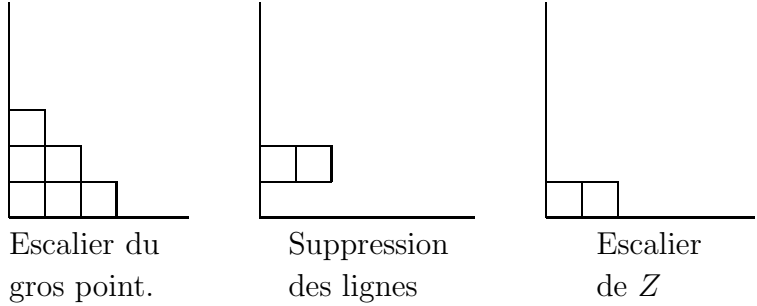
\begin{center}
\setlength{\unitlength}{0.01250000in}%
\begingroup\makeatletter\ifx\SetFigFont\undefined
\def\x#1#2#3#4#5#6#7\relax{\def\x{#1#2#3#4#5#6}}%
\expandafter\x\fmtname xxxxxx\relax \def\y{splain}%
\ifx\x\y   
\gdef\SetFigFont#1#2#3{%
  \ifnum #1<17\tiny\else \ifnum #1<20\small\else
  \ifnum #1<24\normalsize\else \ifnum #1<29\large\else
  \ifnum #1<34\Large\else \ifnum #1<41\LARGE\else
     \huge\fi\fi\fi\fi\fi\fi
  \csname #3\endcsname}%
\else
\gdef\SetFigFont#1#2#3{\begingroup
  \count@#1\relax \ifnum 25<\count@\count@25\fi
  \def\x{\endgroup\@setsize\SetFigFont{#2pt}}%
  \expandafter\x
    \csname \romannumeral\the\count@ pt\expandafter\endcsname
    \csname @\romannumeral\the\count@ pt\endcsname
  \csname #3\endcsname}%
\fi
\fi\endgroup
\begin{picture}(316,127)(13,685)
\thinlines
\put( 15,810){\line( 0,-1){ 90}}
\put( 15,720){\line( 1, 0){ 75}}
\put(120,810){\line( 0,-1){ 90}}
\put(120,720){\line( 1, 0){ 90}}
\put(240,810){\line( 0,-1){ 90}}
\put(240,720){\line( 1, 0){ 87}}
\put( 15,765){\line( 1, 0){ 15}}
\put( 30,765){\line( 0,-1){ 15}}
\put( 30,750){\line( 1, 0){ 15}}
\put( 45,750){\line( 0,-1){ 15}}
\put( 45,735){\line( 1, 0){ 15}}
\put( 60,735){\line( 0,-1){ 15}}
\put( 15,750){\line( 1, 0){ 15}}
\put( 30,750){\line( 0,-1){ 30}}
\put( 15,735){\line( 1, 0){ 30}}
\put( 45,735){\line( 0,-1){ 15}}
\put(120,750){\line( 1, 0){ 30}}
\put(150,750){\line( 0,-1){ 15}}
\put(150,735){\line(-1, 0){ 30}}
\put(240,735){\line( 1, 0){ 30}}
\put(270,735){\line( 0,-1){ 15}}
\put(255,735){\line( 0,-1){ 12}}
\put(135,750){\line( 0,-1){ 15}}
\put(255,726){\line( 0,-1){  6}}
\put( 15,705){\makebox(0,0)[lb]{\smash{\SetFigFont{12}{14.4}{rm}Escalier du }}}
\put( 15,687){\makebox(0,0)[lb]{\smash{\SetFigFont{12}{14.4}{rm}gros point.}}}
\put(132,705){\makebox(0,0)[lb]{\smash{\SetFigFont{12}{14.4}{rm}Suppression }}}
\put(132,687){\makebox(0,0)[lb]{\smash{\SetFigFont{12}{14.4}{rm}des lignes}}}
\put(261,705){\makebox(0,0)[lb]{\smash{\SetFigFont{12}{14.4}{rm}Escalier }}}
\put(261,687){\makebox(0,0)[lb]{\smash{\SetFigFont{12}{14.4}{rm}de $Z$}}}
\end{picture} 
\end{center}
On a vu que l'on avait
besoin successivement d'une puis de trois conditions 
sur la droite $D$. On effectue alors la proc\'edure suivante. 
On supprime dans l'escalier les lignes de longueur un et trois, puis on
``pousse'' les cubes restant vers le bas.
On obtient ainsi un nouvel escalier, associ\'e \`a un 
sous-sch\'ema mon\^omial $Z$ de $\plp$. 
Notons ${\cal C}$ le syst\`eme lin\'eaire form\'e des 
diviseurs $\Delta$ de degr\'e six contenant deux fois $D$ et 
pour lesquels 
$\Delta-2D$ est une courbe de degr\'e quatre 
contenant  $Z$ (contenir $Z$ 
s'interpr\`ete 
g\'eom\'etriquement par le fait que $\Delta-2D$ 
a un ordre de contact d'ordre 
deux avec $D$ en $p_4(0)$). 
Le th\'eor\`eme \ref{theo} 
\'etablit l'in\'egalit\'e $dim \fxl_t \leq dim  \cal C$. 
Puisque $\cal C$ est de dimension douze, on a bien $dim \fxl_t=12$.
\nl
Plus g\'en\'eralement, notre th\'eor\`eme s'int\'eresse \`a la dimension
de certains syst\`emes lin\'eaires $\fxl _t$. On associe \`a $\fxl_t$ un
syst\`eme $\fxc$ au moyen d'op\'erations combinatoires et on \'etablit 
l'in\'egalit\'e $dim \fxl_t \leq dim \fxc$.
\nl
En fait, notre th\'eor\`eme  
ne s'appliquera pas uniquement \`a $\plp$ 
et \`a la sp\'ecia\-li\-sa\-tion de gros 
points sur des droites, comme cela a \'et\'e le cas 
dans l'exemple. Il s'appliquera \`a toute
vari\'et\'e projective irr\'eductible $X$ sur un corps 
al\-g\'e\-bri\-que\-ment clos de caract\'eristique quelconque, et 
\`a la sp\'ecialisation de sch\'emas mon\^omiaux 
en un point $p$
d'un diviseur de Weil irr\'eductible 
$D$ de $X$ en lequel $D$ et $X$ sont lisses.
\nl
La d\'emonstration 
consiste essentiellement \`a contr\^oler  le syst\`eme lin\'eaire
limite (th\'eor\`eme \ref{propfond}) lorsqu'on sp\'ecialise le sch\'ema 
mon\^omial , ce qui est obtenu par une \'etude diff\'erentielle.
\nl
Notre \'enonc\'e est tr\`es similaire \`a la m\'ethode 
d'Horace diff\'erentielle introduite dans  [AH1]. 
D'un c\^ot\'e, le th\'eor\`eme pr\'esent\'e ici 
est plus g\'en\'eral puisque la m\'ethode d'Alexander et Hirschowitz 
ne permet d'utiliser qu'une seule tranche d'un sch\'ema mon\^omial
(i.e. avec les notations du th\'eor\`eme \ref{theo}, ils se limitent 
au cas $r=1$). Mais d'un autre c\^ot\'e, Alexander et Hirschowitz s'autorisent
\`a bouger simultan\'ement plusieurs sch\'emas mon\^omiaux, alors que 
la m\'ethode pr\'esent\'ee  ici ne permet de sp\'ecialiser qu'un unique 
sch\'ema mon\^omial. Avec quelques adaptations dans les 
d\'emonstrations, il aurait \'et\'e 
possible de donner un \'enonc\'e qui englobe l'\'enonc\'e 
d'Alexander-Hirschowitz et le notre.  Mais un tel \'enonc\'e serait 
beaucoup plus  technique et ne donnerait pas lieu \`a de nouvelles 
appplications. 
Sous la forme pr\'esent\'ee dans ce travail,
le th\'eor\`eme \'etait pressenti non seulement par 
Alexander et Hirscho\-witz, mais aussi par 
Simpson qui avait fait une conjecture en ce  sens 
d\`es 1995. 
\nl
Signalons aussi que Joe Harris a donn\'e une conf\'erence 
\`a  Alghero en Juin dernier (1997)
dans laquelle il a annonc\'e  
avoir obtenu avec Lucia Caporaso des 
r\'esultats similaires \`a ceux de cet article
quand $X$ est de dimension deux, 
et quand la caract\'eristique
du corps de base est nulle ou assez grande.
\nl
On trouvera des applications du th\'eor\`eme dans [E] o\`u 
on montre qu'il n'existe pas 
de courbe plane de degr\'e cent soixante quatorze 
 contenant dix points singuliers d'ordre cinquante cinq (ce qui, en un 
 sens \`a pr\'eciser, constitue le premier cas ``critique'' pour lequel la 
postulation de points singuliers ordinaires n'est pas connue).
\nl
Le plan de l'article est le suivant.
Dans la  section \ref{monom}, on explique le lien 
entre les sch\'emas
mon\^omiaux et les escaliers. La section \ref{dechargeable}
est une section technique d'alg\`ebre commutative 
utile pour la d\'emonstration du th\'eor\`eme. 
Le th\'eor\`eme est \'enonc\'e et d\'emontr\'e 
dans la section \ref{thm}.

\section{Sch\'emas mon\^omiaux}
\label{monom}
\subsection{D\'efinition des sch\'emas mon\^omiaux}
On appelle escalier une partie $E$ de $\NN^n$ dont le compl\'ementaire 
$C$ v\'erifie $C + \NN^n \inc C$. Dans la suite, nous ne manipulerons
que des escaliers finis.
On dira par abus de 
langage qu'un mon\^ome $m=x_1^{a_1}x_2^{a_2}\dots x_n^{a_n}$ 
de $k[[x_1,\dots,x_n]]$ 
est dans $E$ si $(a_1,a_2,\dots,a_n)$ est dans $E$. 
L'escalier $E$ d\'efinit un id\'eal
$I^E$ de $k[[x_1,\dots,x_n]]$ qui est  
 l'id\'eal engendr\'e par les mon\^omes hors  de 
$E$. 
\nl
Soit $p$ un point lisse d'une vari\'et\'e $X$ de dimension $n$. 
Le compl\'et\'e $\hat{O_p}$ de l'anneau 
local de $X$ en $p$ est isomorphe \`a l'anneau de  s\'eries formelles 
$k[[x_1,\dots,x_n]]$. Le choix d'un isomorphisme
induit un syst\`eme de coordonn\'ees locales en $p$, not\'e 
$\phi: \s k[[x_1, \dots, x_n]] \fd X$. Moyennant ce choix, 
tout sous-sch\'ema ponctuel de $X$ support\'e par $p$ peut 
\^etre vu comme un sous-sch\'ema de $\s k[[x_1,\dots,x_n]]$. 

\begin{defi}
Un sous-sch\'ema ponctuel $Z$ de $X$ support\'e par $p$ est dit mon\^o\-mial 
d'escalier $E$ si on peut choisir un isomorphisme  
entre $\hat {O_p} $
et $ k[[x_1,\dots,x_n]]$  tel que l'id\'eal d\'efinissant $Z$ dans 
$\s k[[x_1,\dots,x_n]]$ 
soit $I^E$. On notera $X_{\phi}(E)$ le sch\'ema mon\^omial d\'efini par 
$\phi$ et $E$. 
\end{defi}

\bex 
Les gros points de taille $m$ de $X$ sont les sch\'emas mon\^omiaux 
d'escalier $E_m$, avec $E_m=\{(a_1,a_2,\dots,a_n),\ a_1+a_2
+\dots +a_n<m\}$. 
\eex

\subsection{ D\'ecoupage d'un escalier en tranches. Suppression de tranches}
\bdefi \normalfont
Un escalier $E$ de $\NN^n$
d\'efinit une famille d'escaliers $T(E,k)$
de $\NN^{n}$ index\'ee par $\NN-\{0\}$:
 $$T(E,k):=
\{(0,a_2,a_3,\dots,a_n) \mbox{ pour lesquels }
(k-1,a_2,a_3,\dots,a_n) \in E\}$$
L'escalier $T(E,k)$ est  appel\'e $k^{\mbox{\`eme}}$ tranche de 
$E$.
\nl
Un escalier fini peut \^etre caract\'eris\'e par une application 
hauteur $h_E$
de $\NN^{n-1}$ dans $\NN$ qui v\'erifie $h_E(a+b)\leq h_E(a)$ pour tout 
couple $(a,b)$ de $(\NN^{n-1})^2$: l'escalier d\'efini par $h_E$ 
est l'ensemble des $n$-uplets $(a_1,\dots,a_n)$ v\'erifiant 
$a_1 < h_E(a_2,\dots,a_n)$. 
\nl
Pour un escalier $E$ d\'efini par une fonction $h_E$  et un entier $n_i>0$,
on appelle escalier r\'esiduel apr\`es suppression de la  
tranche $n_i$  l'escalier $S(E,n_i)$ d\'efini par la fonction hauteur 
$h_{S(E,n_i)}$:
\beq
h_{S(E,n_i)}(a_2,\dots,a_n)&=&h_E(a_2,\dots,a_n)\mbox{ si }
n_i > h_E(a_2,\dots,a_n)
\\
&=&h_E(a_2,\dots,a_n)-1 \mbox{ si } n_i \leq h_E(a_2,\dots,a_n)
\eeq
\hspace{4cm}
\setlength{\unitlength}{0.00083300in}%
\begingroup\makeatletter\ifx\SetFigFont\undefined%
\gdef\SetFigFont#1#2#3#4#5{%
  \reset@font\fontsize{#1}{#2pt}%
  \fontfamily{#3}\fontseries{#4}\fontshape{#5}%
  \selectfont}%
\fi\endgroup%
\begin{picture}(2274,1856)(6889,-1605)
\thicklines
\put(6901,-699){\line( 1,-1){187.500}}
\put(7089,-886){\line( 0,-1){188}}
\put(7464,-1074){\line( 0, 1){188}}
\put(7464,-886){\line( 1, 1){187}}
\put(7464,-886){\line( 0, 1){187}}
\put(7464,-699){\line( 1, 1){187.500}}
\put(7089,-886){\line( 1, 1){187}}
\put(7276,-699){\line( 0,-1){187}}
\put(7089,-511){\line( 1,-1){375}}
\put(7276,-699){\line( 0, 1){188}}
\put(7276,-511){\line( 1, 1){187.500}}
\put(7089,-324){\line( 1,-1){375}}
\put(7089,-136){\line( 1,-1){187.500}}
\put(7276,-324){\line( 1, 1){188}}
\put(7089, 51){\line( 1,-1){187}}
\put(7276,-136){\line( 1, 1){187.500}}
\put(7276,-136){\line( 0,-1){375}}
\put(7089,-324){\line( 0,-1){187}}
\put(7089,-511){\line( 1,-1){375}}
\put(7464,-886){\line( 0, 1){187}}
\put(7464,-699){\line(-1, 1){375}}
\put(7464,-699){\line( 1, 1){187.500}}
\put(7651,-511){\line( 0,-1){188}}
\put(7651,-699){\line(-1,-1){187}}
\put(7276,-511){\line( 1, 1){187.500}}
\put(7464,-324){\line( 1,-1){187}}
\put(7651,-511){\line(-1,-1){187.500}}
\put(7464,-699){\line(-1, 1){188}}
\put(6901,-699){\line( 1, 1){188}}
\put(7089,-511){\line( 0, 1){562}}
\put(7089, 51){\line( 1, 1){187.500}}
\put(7276,239){\line( 1,-1){188}}
\put(7464, 51){\line( 0,-1){375}}
\put(7464,-324){\line( 1,-1){187}}
\put(7651,-511){\line( 0,-1){375}}
\put(7651,-886){\line(-1,-1){187.500}}
\put(7464,-1074){\line(-1, 1){188}}
\put(7276,-886){\line(-1,-1){187.500}}
\put(7089,-1074){\line(-1, 1){188}}
\put(6901,-886){\line( 0, 1){187}}
\put(7276,-511){\line( 0,-1){188}}
\put(8401,-1589){\makebox(0,0)[lb]{{de la deuxi\`eme tranche}}}
\put(7464,-699){\line( 0,-1){187}}
\put(8401,-699){\line( 1, 1){188}}
\put(8589,-511){\line( 0, 1){375}}
\put(8589,-136){\line( 1, 1){187}}
\put(8776, 51){\line( 1,-1){187.500}}
\put(8964,-136){\line( 0,-1){375}}
\put(8964,-511){\line( 0, 1){  0}}
\put(8964,-511){\line( 1,-1){187.500}}
\put(9151,-699){\line( 0,-1){187}}
\put(9151,-886){\line(-1,-1){187.500}}
\put(8964,-1074){\line(-1, 1){188}}
\put(8776,-886){\line(-1,-1){187.500}}
\put(8589,-1074){\line(-1, 1){188}}
\put(8401,-886){\line( 0, 1){187}}
\put(8401,-699){\line( 1,-1){187.500}}
\put(8589,-886){\line( 0,-1){188}}
\put(8589,-886){\line( 1, 1){187}}
\put(8776,-699){\line( 0,-1){187}}
\put(8776,-699){\line( 1,-1){187.500}}
\put(8964,-886){\line( 1, 1){187}}
\put(8964,-1074){\line( 0, 1){188}}
\put(8589,-511){\line( 1,-1){187.500}}
\put(8776,-699){\line( 1, 1){188}}
\put(8589,-136){\line( 1,-1){187.500}}
\put(8776,-324){\line( 1, 1){188}}
\put(8776,-324){\line( 0,-1){375}}
\put(8589,-324){\line( 1,-1){187}}
\put(8776,-511){\line( 1, 1){187.500}}
\put(8401,-1449){\makebox(0,0)[lb]{Escalier apr\`es suppression}}
\put(6901,-1449){\makebox(0,0)[lb]{Escalier }}
\end{picture}

\nl
 Si 
$(n_1,n_2,\dots,n_r)$ est un $r$-uplet d'entiers v\'erifiant  
$n_1>n_2 \dots >n_r > 0$,
on d\'efinit l'escalier $S(E,n_1,\dots,
n_r)$ obtenu \`a partir de $E$ par suppression des tranches $n_i$ 
r\'ecursivement: $S(E,n_1,\dots,n_r):=S(S(E,n_1\dots,n_{r-1}),n_r)$.
\edefi

\section{Id\'eaux et transporteurs}
\label{dechargeable}
Dans la section pr\'ec\'edente, nous avons d\'efini un id\'eal 
$I^E$ dans $k[[x_1,\dots,x_n]]$, qui correspond
g\'eom\'etriquement \`a un sch\'ema ponctuel. 
Consid\'erons le morphisme de translation $T$:
\beq
T: k[[x_1,\dots,x_n]] &\fd & k[[x_1,\dots,x_n]]\ox k[[t]]\\
x_1&\mapsto& x_1\ox 1 -1 \ox t\\
x_i&\mapsto& x_i\ox 1 \mbox{ si $i>1$}
\eeq
L'id\'eal 
$$J(E):=T(I^E)k[[x_1,\dots,x_n]] \ox k[[t]]$$
d\'efinit une famille plate de sous-sch\'emas 
de $\s k[[x_1,\dots,x_n]]$ param\'etr\'ee par $\s k[[t]]$ qui correspond
g\'eom\'etriquement \`a une translation du sch\'ema ponctuel dans 
la direction $x_1$.
\nl
Au cours de la d\'emonstration du th\'eor\`eme, nous serons grosso
modo amen\'es \`a effectuer les calculs suivants: partant de 
$J_1:=J(E)$, d\'eterminer $J_2:=(J_1:x_1)$, $J_3:=(J_2:x_1)$ \dots
On aimerait en outre que tous les $J_i$ soient de la forme 
$J(F_i)$ pour un escalier $F_i$ de sorte que les id\'eaux 
soient faciles \`a d\'ecrire et \`a manipuler
via leur escalier. Ce n'est 
malheureusement pas le cas. Il est n\'eanmoins possible de donner 
une notion d'id\'eal associ\'e \`a un escalier de sorte 
que tous les id\'eaux soient contr\^ol\'es par 
le fait que ce sont des id\'eaux associ\'es \`a un escalier.
C'est l'objet de la d\'efinition suivante.
\nl
Pour des raisons techniques, nous ne travaillerons pas dans 
$k[[x_1,\dots,x_n]] \ox k[[t]]$, mais dans 
$k[[x_1,\dots,x_n]]/\goth m^s \ox k[[t]]/t^q$ pour diff\'erents 
$s$ et diff\'erents $q$ (o\`u $\goth m$ d\'esigne 
l'id\'eal maximal de $k[[x_1,\dots,x_n]]$ et $s$ et $q$ sont des 
entiers).  Le probl\`eme reste n\'eanmoins le m\^eme, \`a savoir
contr\^oler des calculs de transporteurs \`a l'aide d'escaliers.
\nl
On notera 
\bit
\item 
$r_{qp}^{su}: k[t]/t^{q}\ox k[[x_1,\dots,x_n]]/\goth m^s
\fd k[t]/t^{p}\ox k[[x_1,\dots,x_n]]/\goth m^u$ la projection 
naturelle,
o\`u $p,q,s,u$ sont quatre  
entiers v\'erifiant$0< p\leq q$ et $0< u \leq s$
\item
$J(E,q,s)$ la projection de $J(E)$ dans 
$k[t]/t^{q}\ox k[[x_1,\dots,x_n]]/\goth m^s$
\item 
$I^E$ l'id\'eal de $k[t]/t^{q}\ox k[[x_1,\dots,x_n]]/\goth m^s$ engendr\'e 
par les mon\^omes hors de $E$.
\eit

\bdefi
\label{def ideal d'esc E}
Soient  $q\geq 1$ et $s\geq 1$ deux entiers, et  $E$
un escalier de  $\NN^n$.
Un id\'eal $J$  de $k[t]/t^{q}\ox k[[x_1,\dots,x_n]]/\goth m^s$
est dit id\'eal d'escalier $E$ s'il v\'erifie:
\bit
\item 
$J=I^E$ si $q=1$
\item
si $q>1$ 
\bit
\item
$J\inc I^{T(E,q)}$
\item
pour tout couple $(p,u)$ avec  $0<p<q$, et $0<u <s$,
$r_{qp}^{su}(J:x_1)$ est un id\'eal  
d'escalier  $S(E,q)$
dans $k[t]/t^p\ox k[[x_1,\dots,x_n]]/\goth m^{u}$.
\eit
\eit
\edefi 
Dans cette d\'efinition, les deux premi\`eres propri\'et\'es 
sont les propri\'et\'es vou\-lues pour un id\'eal d'escalier $E$ 
tandis que la troisi\`eme nous assure que la notion est stable 
par calcul de transporteurs.
\nl
L'id\'eal d'escalier $E$ que nous int\'eresse est le suivant:

\bprop
\label{Jd'escE}
\ L'id\'eal $J(E,q,s)$ de $k[[x_1,\dots,x_n]]/\goth m^s \ox 
k[t]/t^q$ 
est un id\'eal d'es\-ca\-lier $E$.
\end{prop}
\noindent
Le reste de la section est  consacr\'e \`a la d\'emonstration de 
cette proposition.
Commen\c cons par le faire  dans le cas $n=1$. 
Notons $E_h$ l'escalier de $\NN$ 
contenant les \'el\'ements inf\'erieurs strictement
\`a $h$. Puisque tout escalier de $\NN$ est de la forme 
$E_h$ pour un certain $h$, la proposition pour $n=1$ dit 
que l'id\'eal $((x_1-t)^h)$ de $k[[x_1]]/x_1^s\ox k[t]/t^q$ 
est un id\'eal d'escalier $E_h$. 
Pour v\'erifier ce fait, la d\'efinition \ref{def ideal d'esc E}
nous invite \`a effectuer des calculs de transporteurs 
et des restrictions.
Lors des
calculs, 
les id\'eaux successifs apparaissant  
se ressemblent 
au sens o\`u ils admettent tous des syst\`emes 
de g\'en\'erateurs similaires. Nous introduisons 
dans la prochaine d\'efinition la notion d'id\'eaux d\'echargeables
de hauteur $H$, qui sont des  id\'eaux admettant 
un ``bon'' syst\`eme 
de g\'en\'erateurs (et bien \'evidemment, 
tous les id\'eaux apparaissant dans les 
calculs sont des id\'eaux d\'echargeables). Et la propri\'et\'e 
fondamentale est que tout id\'eal d\'echargeable de hauteur $H$ 
est un id\'eal d'escalier $E_H$.
\nl
La raison pour laquelle nous avons introduit la notion d'id\'eal 
d'escalier $E$ alors que finalement nous travaillons dans une classe 
d'id\'eaux plus petite, \`a savoir la classe des id\'eaux d\'echargeables
est la suivante: lors de la d\'emonstration du th\'eor\`eme, 
les  propri\'et\'es qui nous int\'eresseront  vraiment 
pour un id\'eal sont celles qui en font un id\'eal d'escalier $E$. On a donc
mis en \'evidence ces propri\'et\'es dans une d\'efinition. 
N\'eanmoins, pour montrer que $J(E,q,s)$ est un id\'eal d'escalier 
$E$, les calculs sont plus 
commodes dans une classe d'id\'eaux 
plus petite (les d\'echargeables) 
dans laquelle les id\'eaux sont contr\^ol\'es par 
un syst\`eme de g\'en\'erateurs.
\nl
En tant que $k$-espace vectoriel, $k[[x_1]]/x_1^s\ox k[t]/t^q$ s'identifie 
au sous-espace vectoriel de $k[x_1,t]$ form\'e par les polyn\^omes 
de degr\'e en $x_1$ plus petit que $s$ et de degr\'e en $t$ plus 
petit que $q$. On dit qu'un \'el\'ement $x_1^{\beta}$ divise 
un \'el\'ement $Q$ de $k[[x_1]]/x_1^s\ox k[t]/t^q$, 
et on \'ecrira $e=\f{Q}{x_1^{\beta}}$
si, moyennant l'identification pr\'ec\'edente, 
$Q$ est une combinaison lin\'eaire 
$\sum \lambda_i x_1^{a_i}t^{b_i}$  de mon\^omes 
o\`u chaque $a_i$ est plus grand que $\beta$, et 
$e=\sum \lambda_i x_1^{a_i-\beta}t^{b_i}$

\bdefi
Un id\'eal $I$ de $k[[x_1]]/x_1^s \ox k[t]/t^q$ est dit d\'echargeable de 
hauteur $H$ s'il est engendr\'e par des \'el\'ements $(e_1,\dots,e_r)$
avec 
\bit
\item
$e_1=\f{(x_1-t)^h}{x_1^{\beta_1}}$ pour des entiers $h$ et $\beta_1$
v\'erifiant $H=h-\beta_1$,  et ($\beta_1=0$ si $q>H$)
\item
pour $i\geq 2$, $e_i=\f{t^{\alpha_i}(x_1-t)^h}{x_1^{\beta_i}}$ 
avec: 
$\alpha_i \geq 1$ et, 
$\forall p\leq q, \ x_1^{q-p+1}$ divise 
$r_{qp}^{ss}(e_i)$.
\eit
\edefi

\bprop
Soit $I=(e_1,\dots,e_r)$ un id\'eal d\'echargeable de hauteur $H$
de $k[[x_1]]/x_1^s \ox k[t]/t^q$.
Si $q\leq H$, alors $(I:x_1)=(x_1^{s-1},\f{e_1}{x_1}, \f{e_2}{x_1},
\dots,\f{e_r}{x_1})$. Si $q>H$, alors 
$(I:x_1)=(x_1^{s-1},e_1,\f{t^{q-h}e_1}{x_1},\f{e_2}{x_1}, \f{e_3}{x_1},
\dots,\f{e_r}{x_1})$.
\end{prop}

\noindent
\textit{D\'emonstration: \\ Le cas $q\leq  H$}: 
$e_1$, vu comme polyn\^ome en $x_1$, admet comme terme constant 
un multiple de $t^H$. Donc ce terme est nul et 
$e_1$ est bien divisible par $x_1$. 
Les \'el\'ements $e_2,\dots,e_r$ sont divisibles par $x_1$ 
par d\'efinition des id\'eaux d\'echargeables.
L'id\'eal $(x_1^{s-1},\f{e_1}{x_1}, \f{e_2}{x_1},
\dots,\f{e_r}{x_1})$ est donc bien d\'efini. 
L'inclusion $(I:x_1) \supset (x_1^{s-1},\f{e_1}{x_1}, \f{e_2}{x_1},
\dots,\f{e_r}{x_1})$ \'etant \'evi\-den\-te, il nous reste \`a voir
qu'un \'el\'ement $m$ de \!$(I:x_1)$ est dans $(x_1^{s-1}\!\!,
\f{e_1}{x_1}, \f{e_2}{x_1},
\dots,\f{e_r}{x_1})$. L'\'el\'ement $x_1m$, qui  est dans $I$, 
s'\'ecrit $\sum \lambda _i e_i$ o\`u les $\lambda_i$ sont  des \'el\'ements 
de $k[[x_1]]/x_1^s \ox k[t]/t^q$. D'o\`u la relation
$$x_1(m-\sum \lambda_i \f{e_i}{x_1})=0$$
Le noyau de la multiplication par $x_1$ \'etant l'id\'eal $(x_1^{s-1})$, 
$m$ est bien dans l'id\'eal 
$(x_1^{s-1},\f{e_1}{x_1},\f{e_2}{x_1}, \f{e_3}{x_1},
\dots,\f{e_r}{x_1})$
\nl
\textit{ Le cas $q>H$}: commme pr\'ec\'edemment, la seule chose non 
imm\'ediate est qu'un \'el\'ement $m$ de $(I:x_1)$ est 
dans l'id\'eal $(x_1^{s-1},e_1,\f{t^{q-h}e_1}{x_1},\f{e_2}{x_1}, 
\f{e_3}{x_1},
\dots,\f{e_r}{x_1})$. Toujours comme pr\'ec\'edemment,
on a l'\'egalit\'e
\beqn
\label{eq7}
x_1.m=\sum \lambda_i {e_i}.
\eeqn
En utilisant l'identification expliqu\'ee plus haut, $\lambda_1$
peut \^etre vu comme un \'el\'ement de $k[x_1,t]$ et on peut \'ecrire 
la division 
$$\lambda_1=x_1.Q+R$$ o\`u $R$ est un \'el\'ement de $k[t]$.
Cette expression et l'expression (\ref{eq7}) fournissent 
l'\'egalit\'e:
$$x_1(m-\sum_{i\geq 2}\lambda_i \f{e_i}{x_1}-Qe_1)=R e_1.$$
Donc $x_1$ divise $R e_1$, ce qui n'est possible que si $R$ est un multiple 
de $t^{q-h}$: $R=\mu t^{q-h}$.
Finalement, l'\'egalit\'e
$$
x_1(m-\sum_{i\geq 2} \lambda_i \f{e_i}{x_1}- 
Qe_1 - {\mu} \f{t^{q-h}e_1}{x_1})=0
$$
et le fait que 
le noyau de la multiplication par $x_1$ est l'id\'eal $(x_1^{s-1})$, 
nous assurent que $m$ est  dans l'id\'eal 
$(x_1^{s-1},e_1,\f{t^{q-h}e_1}{x_1},\f{e_2}{x_1}, \f{e_3}{x_1},
\dots,\f{e_r}{x_1})$.
\findem

\bcor
\label{calcul transp}
Si $I$ est un id\'eal 
d\'echargeable de hauteur $H$ de $k[[x_1]]/x_1^s \ox k[t]/t^q$
et si $q \leq H$, alors pour 
tout couple $(p,u)$ v\'erifiant 
$p<q$ et $u<s$, $r_{qp}^{s u}(I:x_1)$ est un 
id\'eal d\'echargeable de hauteur $H-1$ de  $k[[x_1]]/x_1^u \ox k[t]/t^p$.
\\
Si $I$ est d\'echargeable de hauteur $H$ et si $q > H$, 
$r_{qp}^{su}(I:x_1)$ est un id\'eal d\'echargeable de hauteur $H$
de  $k[[x_1]]/x_1^u \ox k[t]/t^p$.
\ecor
\textit{D\'emonstration:} si $q\leq H$, $r_{qp}^{su}(I:x_1)$ 
admet $(e'_1,\dots,e'_r)$ comme g\'en\'erateurs 
avec $e'_i=r_{qp}^{su}(\f{e_i}{x_1})$. L'\'el\'ement  $e'_1$ v\'erifie 
trivialement la premi\`ere condition demand\'ee aux g\'en\'erateurs 
d'un id\'eal d\'echargeable
de hauteur $H-1$. 
Pour la deuxi\`eme condition, il faut 
voir que pour tout $p'\leq p$ et $i\geq 2$, $x_1^{p-p'+1}$ divise
$r_{pp'}^{uu}\circ r_{qp}^{su}(\f{e_i}{x_1})=r_{qp'}^{su}(\f{e_i}{x_1})$.
Il suffit pour cela de voir que $x_1^{p-p'+2}$ divise 
$r_{qp'}^{ss}(e_i)$.
Or, par hypoth\`ese, $I=(e_1,\dots,e_r)$ est un id\'eal 
d\'echargeable donc  $x_1^{q-p'+1}$ divise $r_{qp'}^{ss}(e_1)$, et 
$q-p'+1\geq p-p'+2$.
\\
Dans le cas $q>H$, $r_{qp}^{su}(I:x_1)$ est de la 
forme $(e'_1,\dots,e'_{r+1})$ avec $e'_1=r_{qp}^{su}(e_1)$,
$e'_i=r_{qp}^{su}(e_i/x_1)$ pour $2\leq i \leq r$ et 
$e'_{r+1}=r_{qp}^{su}(\f{t^{q-h}(x_1-t)^h}{x_1})$. 
Toutes les v\'erifications, sauf une, sont les 
m\^emes qu'au cas pr\'ec\'edent: il nous 
faut en outre montrer que pour tout $p'<p$, $x_1^{p-p'+1}$ divise 
$r_{pp'}^{uu}(e'_{r+1})$. Ceci est vrai car le coefficient en 
$x_1^k$ de $r_{pp'}^{uu}(e'_{r+1})$ est un multiple de $t^{q-1-k}$:
si $k$ est inf\'erieur ou \'egal \`a $p-p'$, 
il  est strictement 
plus petit que $q-p'$, l'exposant $q-1-k$ 
de $t$ est strictement plus grand que $p'-1$ donc $t^{q-1-k}$ est nul
dans $k[[x_1]]/x_1^u \ox k[t]/t^{p'}$.
\findem

\bcor
Si $I=(e_1,\dots,e_r)$ est un id\'eal de 
$k[[x_1]]/x_1^s\ox k[t]/t^q$ d\'echar\-gea\-ble de hauteur $H$, 
alors $I$ est un id\'eal 
d'escalier $E_H$.
\ecor
\textit{D\'emonstration}:
par r\'ecurrence sur $q$. Pour $q=1$, tous les termes $e_i$ avec 
$i\geq 2$ d'un id\'eal d\'echargeable $I=(e_1,\dots,e_r)$ 
sont nuls car ils sont de la forme 
$\f{t^{\alpha_i}(x_1-t)^h}{x_1^{\beta_i}}$ avec $\alpha_i\geq 1$.
Donc $I=(e_1)$ et $e_1=\f{(x_1-t)^h}{x_1^{\beta_1}}=\f{x_1^h}{x_1^{\beta_1}}
=x_1^{H}$. On a bien $I=I^{E_H}$.
\nl
Pour $q>1$, il faut voir que $I$ est inclus dans $I^{T(E_H,q)}$ et 
que, pour $p<q$ et $u<s$,  $r_{qp}^{su}(I:x_1)$ 
est un id\'eal d'escalier $S(E_H,q)$.
\\
Si $q$ est plus grand que $H$, $I^{T(E_H,q)}$ est l'id\'eal unit\'e
donc la premi\`ere condition est trivialement v\'erifi\'ee. 
Dans ce cas, $S(E_H,q)=E_H$. D'apr\`es la proposition \ref{calcul transp}, 
$r_{qp}^{su}(I:x_1)$ est un id\'eal d\'echargeable de hauteur 
$H$, donc c'est un id\'eal d'escalier $E_H$ par hypoth\`ese de 
r\'ecurrence.
\\
Si $q$ est inf\'erieur ou \'egal \`a $H$, la premi\`ere condition dit 
que $I$ est inclus dans l'id\'eal $(x_1)$. V\'erifions que c'est le 
cas pour chacun des g\'en\'erateurs de $I$. C'est vrai pour les 
\'el\'ements $e_2,\dots,e_r$ par d\'efinition des g\'en\'erateurs 
d'un id\'eal d\'echargeable. C'est \'egalement vrai pour $e_1=
\f {(x_1-t)^h}{x_1^{\beta_1}}$ car son terme constant est un 
multiple de $t^H$, donc est nul. 
\\
Pour la deuxi\`eme condition, il faut voir que 
$r_{qp}^{su}(I:x_1)$  est un id\'eal d'escalier 
$S(E_H,q)=E_{H-1}$. Or, d'apr\`es  
la proposition \ref{calcul transp}, $r_{qp}^{su}(I:x_1)$ 
 est un id\'eal d\'echargeable de hauteur $H-1$. C'est donc aussi un 
 id\'eal d'escalier $E_{H-1}$ par 
l'hypoth\`ese de 
r\'ecurrence.
\findem

\bcor
\label{J d'esc E cas n=1}
Soit $E_h \inc \NN$ un escalier. L'id\'eal $J(E_h,s,q)$ 
de $k[[x_1]]/x_1^s\ox k[t]/t^q$ est un id\'eal 
d'escalier $E_h$.
\ecor
\textit{D\'emonstration}: $J(E_h,s,q)=((x_1-t)^h)$ est trivialement 
un id\'eal d\'echar\-gea\-ble de hauteur $h$. 
C'est donc un id\'eal d'escalier 
$E_h$ d'apr\`es le corollaire pr\'ec\'edent.
\findem
Soit $E$ un escalier de $\NN^n$. 
On va maintenant montrer pour $n$ quelconque  
que $J(E,s,q)$ est un id\'eal 
d'escalier $E$ de $k[[x_1,\dots,x_n]]/\goth m^s \ox k[t]/t^q$ 
en se 
ramenant au cas  $n=1$. Identifions pour cela ensemblistement
l'anneau
$k[[x_1,\dots,x_n]]/\goth m^s \ox k[t]/t^q$ au produit 
$$
 \prod_{(\alpha_2,\dots,\alpha_n)\; t.q.\ 
s-\alpha_2-\dots-\alpha_n\geq 1}k[[x_1]]/x_1^{s-\alpha_2-\dots-\alpha_n}\ox k[t]/t^q
$$
o\`u l'identification envoie un terme $m$ de la composante 
d'indice $(\alpha_2,\dots,\alpha_n)$ sur le produit 
$mx_2^{\alpha_2}x_3^{\alpha_3}\dots x_n^{\alpha_n}$.

\blm
\label{J est gradue}
Soient $E$ un escalier fini de $\NN^n$ d\'efini par une fonction 
hauteur $h_E$ 
et $I_{\alpha_2,\dots,\alpha_n}$ l'id\'eal de 
$k[[x_1]]/x_1^{s-\alpha_2-\dots-\alpha_n}\ox k[t]/t^q$ engendr\'e par 
$(x_1-t)^{h_E(\alpha_2,\dots,\alpha_n)}$.  
L'id\'eal $J(E,s,q)$ co\"\i ncide avec le produit
$\prod_{\alpha_2,\dots,\alpha_n}I_{\alpha_2,\dots,\alpha_n}$
\elm
\textit{D\'emonstration}: puisque chaque $I_{\alpha_2,\dots,\alpha_n}$
est inclus dans $J(E,s,q)$, on a l'inclusion 
$$ \prod_{\alpha_2,\dots,\alpha_n}I_{\alpha_2,\dots,\alpha_n}
\inc J(E,s,q)$$ 
Les \'el\'ements 
$(x_1-t)^{h_E(\alpha_2,\dots,\alpha_n)}
x_2^{\alpha_2}.x_3^{\alpha_3}.\dots.x_n^{\alpha_n} $ 
engendrent $J(E,s,q)$ et sont dans
$ \prod_{\alpha_2,\dots,\alpha_n}I_{\alpha_2,\dots,\alpha_n}$ . 
Il suffit donc  pour montrer 
l'inclusion inverse de v\'erifier que le produit 
$\prod_{\alpha_2,\dots,\alpha_n}I_{\alpha_2,\dots,\alpha_n}$
est un id\'eal de $k[[x_1,\dots,x_n]]/\goth m^s \ox k[t]/t^q$. 
Ce produit est clairement un 
$k[t]/t^q$-module. Utilisant alors 
la lin\'earit\'e, il suffit de v\'erifier que
le produit d'un \'el\'ement $e_0$ de $I_{\alpha_2^0,\dots,\alpha_n^0}$
et d'un mon\^ome 
$m=x_1^{\beta_1}.x_2^{\beta_2}.x_3^{\beta_3}.\dots.x_n^{\beta_n}$
est dans $\prod_{\alpha_2,\dots,\alpha_n}I_{\alpha_2,\dots,\alpha_n}$.
Par d\'efinition de $I_{\alpha_2^0,\dots,\alpha_n^0}$, 
$$e_0=(x_1-t)^{h_E(\alpha_2^0,\dots,\alpha_n^0)}.
x_2^{\alpha_2^0}.x_3^{\alpha_3^0}.\dots.x_n^{\alpha_n^0}.\mu$$
o\`u $\mu$ est 
un \'el\'ement de $k[[x_1]]/x_1^s\ox k[t]/t^q$.
On a donc $$m.e_0= 
(x_1-t)^{h_E(\alpha_2^0,\dots,\alpha_n^0)+\beta_1}.
x_2^{\alpha_2^0+\beta_2}.x_3^{\alpha_3^0+\beta_3}.\dots.
x_n^{\alpha_n^0+\beta_n}.\mu$$ 
Le terme $me_0$ est donc aussi 
un multiple de 
$$(x_1-t)^{h_E(\alpha_2^0+\beta_2,\dots,\alpha_n^0+\beta_n)}.
x_2^{\alpha_2^0+\beta_2}.x_3^{\alpha_3^0+\beta_3}.\dots.
x_n^{\alpha_n^0+\beta_n}$$ 
en vertu de l'in\'egalit\'e $$h_E(\alpha_2^0+\beta_2,
\alpha_3^0+\beta_3,\dots,
\alpha_n^0+\beta_n)\leq  
h_E(\alpha_2^0,
\alpha_3^0,\dots,
\alpha_n^0).$$ Par suite $m.e_0$ est dans $I_{\alpha_2^0+\beta_2, \dots,
\alpha_n^0+\beta_n}$.
\findem

\blm
\label{etre d'esc E est une prop graduee}
Soit $E$ un escalier de $\NN^n$
donn\'e par une fonction hauteur $h_E$.
Soit $J$ un id\'eal de $k[[x_1,\dots,x_n]]
/\goth m^s \ox k[t]/t^q$  tel que $J=
\prod_{\alpha_2,\dots,\alpha_n}J_{\alpha_2,\dots,\alpha_n}$,
o\`u chaque $J_{\alpha_2,\dots,\alpha_n}$ est un id\'eal 
de 
$k[[x_1]]
/x_1^{s-\alpha _2-\dots-\alpha_n} \ox k[t]/t^q$
d'escalier  $E_{h_E(\alpha_2,\dots,\alpha_n)}$ . Alors $J$ est un id\'eal
d'escalier $E$.
\elm
\textit{D\'emonstration}: appelons id\'eal gradu\'e de $k[[x_1,\dots,x_n]]
/\goth m^s \ox k[t]/t^q$ un id\'eal $K$ qui s'\'ecrit comme produit 
d'id\'eaux $K=\prod K_{\alpha_2, \dots,\alpha_n}$. On dira que les 
$K_{\alpha_2,\dots,\alpha_n}$ sont les parties gradu\'ees
de $K$. Deux id\'eaux 
gradu\'es $L$ et $K$ v\'erifient $L\inc K$ si et seulement si 
pour tout $(\alpha_2,\dots,\alpha_n)$, $L_{\alpha_2,\dots,\alpha_n}
\inc K_{\alpha_2,\dots,\alpha_n}$. De plus, si $K$ est gradu\'e,
les id\'eaux $(K:x_1)$ et $r_{qp}^{su}(K)$ sont gradu\'es et, plus
pr\'ecis\'ement, 
$(K:x_1)=\prod (K_{\alpha_2,\dots,\alpha_n}:x_1)$ et 
$r_{qp}^{su}(K)=\prod r_{qp}^{su}(K_{\alpha_2,\dots,\alpha_n})$. En 
d\'efinitive, dans la d\'efinition \ref{def ideal d'esc E}, 
toutes les v\'erifications 
\`a faire concernent des id\'eaux gradu\'es, et les calculs de 
transporteur et les restrictions respectent la graduation.  
Donc le fait d'\^etre 
un id\'eal d'escalier $E$ se v\'erifie sur chaque partie 
gradu\'ee.
\findem

\noindent
\textit{D\'emonstration de la proposition \ref{Jd'escE}}:
d'apr\`es le lemme \ref{J est gradue}, l'id\'eal $J(E,s,q)$ est un produit 
d'id\'eaux $I_{\alpha_2,\dots,\alpha_n}$. Chacun de ces id\'eaux  
$I_{\alpha_2,\dots,\alpha_n}$ est un id\'eal d'escalier 
$E_ {h(\alpha_2,\dots,\alpha_n)}$ d'apr\`es le corollaire
 \ref{J d'esc E cas n=1}. 
On conclut enfin avec le lemme 
\ref{etre d'esc E est une prop graduee} que 
$J(E,s,q)$ est un id\'eal d'escalier $E$.
\findem

\section{Le th\'eor\`eme}
\label{thm}
Le th\'eor\`eme traite de syst\`emes lin\'eaires. Comme 
dans l'exemple introductif, les syst\`emes  consid\'er\'es
seront des sous-syst\`emes $\fxl_t$
d'un syst\`eme lin\'eaire $\fxl$; les diviseurs de $\fxl_t$ seront 
des diviseurs de $\fxl$ qui contiennent un sch\'ema mon\^omial 
$X(t)$ variant avec le temps $t$.  Au temps $t=0$, le sch\'ema 
$X(0)$ se sp\'ecialise sur un diviseur de Weil irr\'eductible $D$.
Le trajet du sch\'ema mon\^omial 
sera une translation relativement \`a un syst\`eme de coordonn\'ees 
locales ``compatible'' avec le diviseur $D$. Expliquons 
ce que cela signifie. 
\nl
Soient $X$ une vari\'et\'e projective irr\'eductible
de dimension $n$, $D$ une sous-vari\'et\'e irr\'eductible 
de $X$ de dimension $n-1$. 
Soit $p$ un point de $D$ en lequel  $X$ et $D$ sont lisses. 
Choisissons une fois pour toutes un syst\`eme de coordonn\'ees locales 
$\phi:\s k[[x_1,\dots,x_n]] \fd X $ en $p$ de sorte que $D$ soit localement 
d\'efini par $x_1=0$. 
L'id\'eal 
$J(E)$
introduit au d\'ebut de la section pr\'ec\'edente
d\'efinit une famille plate de sous-sch\'emas 
de $\s k[[x_1,\dots,x_n]]$ param\'etr\'ee par $\s k[[t]]$. On peut 
\'egalement voir cette famille comme une famille plate de sous-sch\'emas 
de $X$ moyennnant le morphisme de coordonn\'ees locales  $\phi$. 
On note $X_{\phi}(E,t)$ 
la fibre g\'en\'erique de cette famille plate. La fibre sp\'eciale 
de cette famille est $X_{\phi}(E,0)=X_{\phi}(E)$. Cette famille plate est
associ\'ee \`a un morphisme $\s k[[t]]\fd Hilb(X)$ qui 
correspond au trajet du sch\'ema mon\^omial d\'efini par la translation.
\nl
Soient $\fxl$ un syst\`eme lin\'eaire de diviseurs de Cartier sur $X$
et $Y$ un sous-sch\'ema de $X$. On note
$\fxl(-Y)$ le sous-syst\`eme lin\'eaire de $\fxl$ form\'e par les 
diviseurs qui contiennent $Y$. Si $Y$ et $Z$ sont deux sous-sch\'emas 
de $X$, le produit des id\'eaux $I(Y)$ et $I(Z)$ de $Y$ et $Z$ 
d\'efinit un sous-sch\'ema de $X$ not\'e $Y+Z$. En particulier,
$\fxl(-Y-Z)$ est bien d\'efini, m\^eme si $Y$ est un diviseur de $X$
et $Z$ un sous-sch\'ema de dimension z\'ero.
\nl
Notons $Z_k$ le sous-sch\'ema de $X$ d\'efini par 
le syst\`eme  de coordonn\'ees locales 
$\phi$ et la tranche $T(E,k)$: 
$Z_k:=X_{\phi}(T(E,k))$. 
Les sch\'emas $Z_k$ sont inclus dans le diviseur $D$. 
\\
Les sch\'emas mon\^omiaux d'escalier $E$
s'organisent en une vari\'et\'e irr\'eductible [H]
et on peut donc parler du  sch\'ema g\'en\'erique d'escalier $E$,
qu'on note $X(E)$.

\bthm
\label{theo}
Soient $\fxl$ un syst\`eme lin\'eaire 
sur 
$X$ et  $n_1,n_2,\dots,n_r$ des entiers 
v\'erifiant $n_1>n_2>\dots>n_r>0$. Supposons que  
pour tout $i$ compris entre un et $r$, 
$\fxl(-(i-1) D-Z_{n_i})=\fxl(-iD)$. 
Alors 
$$dim \; \fxl(-X(E)) \leq dim\; \fxl(-rD-X_{\phi}
(\;S(E,n_1,\dots, n_r)\;))$$
\end{thm} \noindent
Le th\'eor\`eme est une cons\'equence imm\'ediate du th\'eor\`eme 
suivant:

\bthm
\label{propfond}
Soient $\fxl$ un syst\`eme lin\'eaire 
sur 
$X$ et  $n_1,n_2,\dots,n_r$ des entiers 
v\'erifiant $n_1>n_2>\dots>n_r>0$. Supposons que  
pour tout $i$ compris entre un et $r$, 
$\fxl(-(i-1) D-Z_{n_i})=\fxl(-iD)$. 
Alors on a l'inclusion
$$\lim_{t \fd 0} \fxl(-X_{\phi}(E,t)) \inc \fxl(-rD-X_{\phi}
(\;S(E,n_1,\dots, n_r)\;))$$
\end{thm}
\noindent
\textit{D\'emonstration du th\'eor\`eme \ref{theo}}: puisque 
$X_{\phi}(E,t)$ est une sp\'ecialisation de $X(E)$, on a par semi-continuit\'e
$$dim\;\fxl(-X(E))\leq dim\; \fxl(-X_{\phi}(E,t))
$$
La limite \'etant par d\'efinition une limite dans une 
Grassmannienne, on a:
$$
dim\;\fxl(-X_{\phi}(E,t)) = dim\;lim_{t \fd 0}\fxl(-X_{\phi}(E,t))
$$
Enfin, la proposition \ref{propfond} implique:
$$
dim\; lim_{t \fd 0}\fxl(-X_{\phi}(E,t)) \leq dim\; \fxl(-rD-X_{\phi}
(\;S(E,n_1,\dots, n_r)\;))
$$ 
Ces in\'egalit\'es mises bout \`a bout donnent l'in\'egalit\'e du 
th\'eor\`eme.
\findem
\textit{D\'emonstration du th\'eor\`eme  \ref{propfond}}:
\nl
Le syst\`eme lin\'eaire $\fxl$ est de la forme $\PP(V)$ pour un 
fibr\'e en droites $F$ sur $X$ et un espace vectoriel $V$  de sections 
de $F$. Notons
$n-1$ la dimension projective 
du syst\`eme lin\'eaire $\fxl(-X_{\phi}(E,t))$.
Il existe un unique morphisme 
$$f: \s k[[t]] \fd \GG(n,V)$$ 
qui envoie le point g\'en\'erique 
sur le point (non ferm\'e) de la grassmannienne 
param\'etrant le 
syst\`eme lin\'eaire $\fxl(-X_{\phi}(E,t))$.
L'image du point sp\'ecial d\'efinit un sous-espace vectoriel $W$ 
de $V$ et, par d\'efinition, $\PP(W)=\lim_{t \fd 0} \fxl(-X_{\phi}(E,t))$
\nl
Restreignons la base du  fibr\'e $F$ \`a $\s \hat O_{X,p}$, o\`u 
$\hat O_{X,p}$ est le compl\'et\'e de l'anneau local de $X$ en 
$p$. Au dessus de cette base, le faisceau localement libre 
$F$ est trivial et on peut 
en choisir un g\'en\'erateur local $g$. Une fois $g$ 
choisi, on peut r\'ealiser toute section 
de $F$ comme une fonction de $\hat O_{X,p}$. 
Le syst\`eme $\phi$ de coordonn\'ees locales en $p$ \'etant donn\'e,
toute fonction de $\hat O_{X,p}$ s'identifie \`a un 
\'el\'ement de $k[[x_1,\dots,x_n]]$.
On dispose donc d'un morphisme, injectif
car $X$ est irr\'eductible:
$$i:V \fd k[[x_1,\dots,x_n]]$$
Notons $p_s$ la projection de $k[[x_1,\dots,x_n]]$ dans $k[[x_1,\dots,x_n]]/
\goth m^s$. 
Puisque $V$ est de dimension finie, le morphisme 
$$p_s \circ i: V \fd k[[x_1,\dots,x_n]]/\goth m ^s$$
est \'egalement injectif pour $s$ assez grand.
Un \'el\'ement $f$ de $V$ s'annule $n$ fois sur $D$ si et seulement si 
$i(f)$ est divisible par $x_1^n$.
Toujours pour $s$ assez grand, 
$f$ s'annule $n$ fois sur $D$ si et seulement si $p_s\circ i(f)$ 
est un multiple de  $x_1^n$. 
\nl
Pour $q \geq 0$, 
notons $f_q$ la restriction du morphisme $f$ \`a  $\s k[t]/t^{q}$: 
$$f_q: \s k[t]/t^{q} \fd \GG(n,V)$$
L'image inverse par 
$$f_q\x Id:\s k[t]/t^{q} \x V \fd \GG(n,V)\x V$$ 
du fibr\'e universel au dessus de 
$\GG(n,V)$ est un sous-fibr\'e $F_q$ de rang $n$ de $\s k[t]/t^{q} \x V$.
Expliquons comment associer un id\'eal 
$I(s_q,s)$ de $k[[x_1,\dots,x_n]]/\goth m^s\ox k[t]/t^{q}$
\`a une section $s_q$ de $F_q$.
\\
Toute section $s_q$ de $F_q$ est  aussi
une section de $\s k[t]/t^{q} \x V$, 
et  est  d\'efinie par  un morphisme de $\s k[t]/t^{q}$ dans $V$.
Par composition \`a droite avec le morphisme $p_s \circ i$, la section 
$s_q$ d\'efinit un morphisme $f(s_q,s)$:
 $$f(s_q,s):\s k[t]/t^{q} \fd k[[x_1,\dots,x_n]]/ \goth m^s.$$
Il existe un ferm\'e  $U$ de $k[[x_1,\dots,x_n]]/ \goth m^s \x 
\s k[[x_1,\dots,x_n]]/ \goth m^s$ dont la fibre au dessus d'un point 
$f$ est le sous-sch\'ema de $\s k[[x_1,\dots,x_n]]/ \goth m^s$ d\'efini 
par l'id\'eal $(f)$. L'image inverse de $U$ par 
$f(s_q,s)\x Id$
est 
un ferm\'e $U(s_q,s)$ de $\s k[t]/t^{q} \x \s k[[x_1,\dots,x_n]]/ \goth m^s$. 
On note $I(s_q,s)$ l'id\'eal de $k[t]/t^{q} \ox k[[x_1,\dots,x_n]]/ \goth m^s$
d\'efinissant $U(s_q,s)$.
\nl
La signification g\'eom\'etrique de $I(s_q,s)$ est la suivante.
La section $s_q$ d\'efinit une famille de diviseurs de $X$ param\'etr\'ee 
par $\s k[t]/t^{q}$, donc un sous-sch\'ema $Z$ de $\s k[t]/t^{q} \x X$.
La trace de $Z$ sur 
$$\s k[t]/t^{q} \x \s k[[x_1,\dots,x_n]]/\goth m^s$$
est un sous-sch\'ema d\'efini par l'id\'eal $I(s_q,s)$.
\nl
Le comportement par restriction des id\'eaux $I(s_q,s)$ est agr\'eable:
si $p,q,s,u$ sont quatre entiers avec $q\geq p$, $s\geq u$, 
et si $s_p$ est la 
restriction de $s_q$ au dessus de  
$\s k[t]/t^{p}$, alors $I(s_p,u)=r_{qp}
^{su}(I(s_q,s))$.
\nl
Pour montrer le th\'eor\`eme, il nous faut voir $(**)$ que pour toute section 
$s_1$ de $F_1$ 
au dessus du point ferm\'e et pour tout entier $s$ assez grand, 
$I(s_1,s)\inc x_1^r. I^{S(E,n_1,\dots,n_r)}$.
\nl
Toute section $s_1$ de $F_1$ 
au dessus du point ferm\'e est la restriction
d'une section $s_{n_1}$ de $F_{n_1}$ 
au dessus de $\s k[t]/t^{n_1}$. Notons 
$s_{n_i}$ la restriction de $s_{n_1}$ \`a $\s k[t]/t^{n_i}$. 
\nl
Montrons la proposition $(*)$
suivante, qui impliquera facilement $(**)$ et donc le th\'eor\`eme \ref{propfond}:
pour $s$ assez grand, 
l'id\'eal 
$I(s_{n_i},s)$ est inclus dans un id\'eal $x_1^i. M(n_i,s)$, o\`u 
pour tout $p<n_i$ et $u<s$, $r_{n_ip}^{su}M(n_i,s)$ 
est un id\'eal 
d'escalier $S(E,n_1,\dots,n_i)$.
\\
On proc\`ede  par r\'ecurrence sur $i$.
\nl
Pour $i=1$, 
on peut dire informellement que $s_{n_1}$ 
est une famille de sections de $F$ param\'etr\'ee par un temps 
$t$ dans $\s k[t]/t^{n_1}$, et que cette famille de 
sections s'annule ``\`a tout instant $t$ sur 
le translat\'e par $t$ dans la direction $x_1$ du sch\'ema 
$X_{\phi}(E)$''. Plus rigoureusement, on a l'inclusion
\beqn
\label{eq1}
I(s_{n_1},s) \inc J(E,n_1,s)
\eeqn
De plus,  $J(E,n_1,s)$ est un id\'eal 
d'escalier $E$ de $k[[x_1,\dots,x_n]]/\goth m^s \ox k[t]/t^q$
d'apr\`es la proposition \ref{Jd'escE}, donc  
\beqn
\label{eq2}
J(E,n_1,s)\inc I^{Tr(E,n_1)}
\eeqn
Pour $s$ assez grand, les inclusions 
(\ref{eq1}) et (\ref{eq2}) montrent que $s_{n_1}$ d\'efinit 
une famille de sections de $F$ s'annulant \`a tout instant
sur $Z_{n_1}$. Donc, par hypoth\`ese, c'est \'egalement une famille de 
sections s'annulant sur $D$. Remarquons qu'\`a priori, 
l'hypoth\`ese dit qu'une section de $F$ qui s'annule sur  $Z_{n_1}$ 
s'annule sur $D$ mais ne dit rien pour les familles de sections.
Cependant, si on note $W_{n_1}$ le lieu sch\'ematique 
dans $V$ form\'e par les sections de $F$ qui s'annulent sur 
$Z_{n_1}$ et $W_D$ le lieu sch\'ematique form\'e par les 
sections qui s'annulent sur $D$, $W_{n_1}$ et $W_D$ sont 
r\'eduits car ce sont des espaces vectoriels. En particulier, 
l'inclusion ensembliste de $W_{n_1}$ dans $W_D$, v\'erifi\'ee 
par hypoth\`ese, 
implique l'inclusion sch\'ematique. Les familles de sections 
de $F$ param\'etr\'ees par une base $B$ et 
s'annulant sur $Z_{n_1}$ correspondent aux morphismes de $B$
dans $W_{n_1}$, qui sont aussi des morphismes de $B$ 
dans $W_D$. Les familles 
de sections s'annulant sur $Z_{n_1}$ s'annulent donc sur 
$D$.
\nl
Puisque $s_{n_1}$ est une famille de sections de $F$ s'annulant 
sur $D$,  tout \'el\'ement $e$ de $I(s_{n_1},s)$ est 
divisible par $x_1$:
$e=x_1.f$, et d'apr\`es la relation \ref{eq1}, 
$f \in (J(E,n_1,s):x_1)$, ce qui s'\'ecrit aussi
$$
I(s_{n_1},s) \inc x_1.(J(E,n_1,s):x_1)
$$
Posons  $M(n_1,s):=(J(E,n_1,s):x_1)$.
Puisque $J(E,n_1,s)$ est un id\'eal d'escalier $E$ et 
par d\'efinition des id\'eaux d'escaliers $E$, pour tout 
$p<n_1$ et $u<s$, 
$r_{n_1p}^{su}(M(n_1,s))=r_{n_1p}^{su}(J(E,n_1,s):x_1)$ 
est bien un id\'eal 
d'escalier $S(E,n_1)$. La proposition $(*)$ est vraie pour $i=1$. 
\nl
Supposons $(*)$ vraie au rang $q-1$. L'inclusion 
\beq
I(s_{n_{q-1}},s+q)\inc x_1^{q-1}.M(n_{q-1},s+q)
\eeq
est v\'erifi\'ee pour $s$ assez grand et implique par la 
restriction $r_{n_{q-1}n_{q}}^{(s+q)(s+q-1)}$
\beqn
\label {eq22}
I(s_{n_q},s+q-1)   \inc x_1^{q-1}.
r_{n_{q-1}n_{q}}^{(s+q)(s+q-1)}M(n_{q-1},s+q)
\eeqn
Puisque $r_{n_{q-1}n_{q}}^{(s+q)(s+q-1)}M(n_{q-1},s+q)$ est un 
id\'eal d'escalier $S(E,n_1,\dots,n_{q-1})$
de $k[[x_1,\dots,x_n]]/\goth m^s \ox k[t]/t^{n_q}$, il est inclus 
dans $I^{T(S(E,n_1,\dots,n_{q-1}),n_q)}=I^{T(E,n_q)}$.
Pour $s$ assez grand, 
l'in\-clu\-si\-on (\ref{eq22})
montre alors que $s_{n_q}$ est une famille 
de sections de $F$ s'an\-nu\-lant sur \!$(q-1)D+Z_{n_q}$, 
ce qui par hypoth\`ese
est aussi une famille de sections de $F$ s'annulant sur $qD$. 
Tout \'el\'ement $e$ de $I(s_{n_q},s+q-1)$ est donc un multiple
de $x_1^q$:
\beqn
\label{eq222}
e=x_1^q.f.
\eeqn
Par l'inclusion (\ref{eq22}), $e$ s'\'ecrit aussi
$$
e=x_1^{q-1}.g
$$ o\`u 
$g$ est dans $r_{n_{q-1}n_{q}}^{(s+q)(s+q-1)}M(n_{q-1},s+q)$.
On en d\'eduit l'\'egalit\'e 
$$x_1^{q-1}(g-x_1f)=0
$$
Puisque le 
noyau de la multiplication par $x_1^{q-1}$ dans 
$k[[x_1,\dots,x_n]]/ \goth m^{q+s-1}\ox k[t]/t^{n_q}$
est inclus dans 
$\goth m^{s}$, 
on a donc 
$$
r_{n_q,n_q}^{q+s-1,s}(f).x_1=r_{n_q,n_q}^{q+s-1,s}(g).
$$
Le terme 
$r_{n_q,n_q}^{q+s-1,s}(g)$ est dans 
$$r_{n_q,n_q}^{q+s-1,s}
 \circ r_{n_{q-1}n_{q}}^{(s+q)(s+q-1)}M(n_{q-1},s+q)=
r_{n_{q-1}n_q}^{s+q,s}M(n_{q-1},s+q)
$$
d'o\`u
$$
r_{n_q,n_q}^{q+s-1,s}(f)\in 
(r_{n_{q-1}n_q}^{s+q,s}M(n_{q-1},s+q):x_1)
$$
L'image de l'\'egalit\'e (\ref{eq222}) par $r_{n_q,n_q}^{q+s-1,s}$
montre alors que 
$$
I(s_{n_q},s)\inc x_1^q.(r_{n_{q-1}n_q}^{s+q,s}M(n_{q-1},s+q)\;:\;x_1).
$$
Posons $M(n_q,s)=(\ 
 r_{n_{q-1}n_q}^{s+q,s}M(n_{q-1},s+q):x_1\ )
$.
On a bien 
$I(s_{n_q},s)$ qui est inclus dans l'id\'eal $x_1^q. M(n_q,s)$. Il 
reste \`a voir que pour 
tout $p<n_q$ et $u<s$, $r_{n_qp}^{su}M(n_q,s)$ 
est un id\'eal 
d'escalier $S(E,n_1,\dots,n_q)$. Ce qui est vrai 
car  $r_{n_qp}^{su}M(n_q,s)
= r_{n_qp}^{su}(\ 
 r_{n_{q-1}n_q}^{s+q,s}M(n_{q-1},s+q):x_1\ )$ 
 et $ r_{n_{q-1}n_q}^{s+q,s}M(n_{q-1},s+q)$ 
 est un id\'eal d'escalier $S(E,n_1,\dots,n_{q-1})$ de $k[[x_1,\dots,x_n]]
 \ox k[t]/t^{n_q}$. 
\nl
La d\'emonstration de la r\'ecurrence est termin\'ee. D\'eduisons 
maintenant $(**)$ de la proposition $(*)$, ce qui ach\`evera la
d\'emonstration du th\'eor\`eme \ref{propfond}.
\nl
Si $n_r\neq 1$, la proposition $(*)$
appliqu\'ee \`a $i=r$ dit que $I(s_{n_r},s+1)$ est inclus dans un 
id\'eal produit $x_1^r.M(n_r,s+1)$. L'image de cette inclusion 
par l'application de restriction $r_{n_r1}^{(s+1)s}$
donne  $I(s_1,s)\inc x_1^r. I^{S(E,n_1,
\dots,n_r)}$ car la restriction de $M(n_r,s+1)$  
est  $I^{S(E,n_1,\dots,n_r)}$
par d\'efinition de $M(n_r,s+1)$ et des id\'eaux  d'escalier 
$ S(E,n_1,\dots,n_r)$. La proposition $(**)$ est donc 
d\'emontr\'ee pour  $n_r\neq 1$.
\nl
Si $n_r=1$, la proposition $(*)$ appliqu\'ee \`a $i={r-1}$ 
dit que $I(s_{n_{r-1}},s+r+1)$ est inclus dans un 
id\'eal $x_1^{r-1}.M(n_{r-1},s+r+1)$. L'image de cette inclusion 
par l'application de restriction $r_{n_{r-1}1}^{(s+r+1)(s+r)}$
est 
\beqn
\label{eq33}
I(s_1,s+r)\inc x_1^{r-1}. I^{S(E,n_1,
\dots,n_{r-1})}.
\eeqn
Puisque $I^{S(E,n_1,
\dots,n_{r-1})}\inc I^{T(E,1)}$, on a \'egalement 
l'inclusion $I(s_1,s+r)\inc x_1^{r-1}.I^{T(E,1)}$,
ce qui signifie pour $s$ assez grand 
que $s_1$ est section de $F$ qui s'annule
sur $(r-1)D+Z_1$.  Par hypoth\`ese, $s_1$
est alors une  section de $F$ qui s'annule sur $rD$. 
Tout \'el\'ement $e$ de $I(s_{1},s+r)$ est donc divisible par $x_1^r$:
\beqn
\label{eq333}
e=x_1^r.f.
\eeqn 
Cette \'egalit\'e, la relation (\ref{eq33}), 
et le fait que le noyau de la multiplication 
par $x_1^{r-1}$ dans $k[[x_1,\dots,x_n]]/\goth m^{s+r}$
soit inclus dans $\goth m^{s+1}$ montrent 
que 
$$
r_{11}^{s+r,s+1}(f) \in (I^{S(E,n_1,
\dots,n_{r-1})}:x_1)=I^{S(E,n_1,
\dots,n_r)}+(x_1^s).
$$
Cette appartenance et la relation (\ref{eq333}) donnent finalement 
$$I(s_1,s) \inc x_1^r.I^{S(E,n_1,
\dots,n_r)}.$$
La d\'emonstration de $(**)$ est termin\'ee. 
\findem
Bibliographie:\nl
[AH1]: Alexander J. et Hirschowitz A., An asymptotic vanishing theorem for 
generic unions of multiple points, duke e-print 9703037
\nl
[AH2]: Alexander J. et Hirschowitz A., La m\'ethode d'Horace \'eclat\'ee: 
application \`a l'interpolation en degr\'e quatre, Invent. Math. 107, (1992),
586-602
\nl
[CM]: Ciliberto C. et Miranda R., On the dimension of linear systems of 
plane curves with general multiple base points, duke e-print 9702015
\nl
[E]: Evain L., Une g\'en\'eralisation de la conjecture de 
Harbourne-Hirschowitz aux points infiniment voisins, pr\'eprint en 
pr\'eparation.
\nl
[H]: Hirschowitz A., La m\'ethode d'Horace pour l'interpolation \`a plusieurs
variables, Manuscripta math, vol. 50, (1995), 337-388

\end{document}